\documentclass[letterpaper, 10 pt, conference]{ieeeconf}  

\IEEEoverridecommandlockouts                              
\usepackage{graphics} 
\usepackage{epsfig} 
\usepackage{mathptmx} 
\usepackage{times} 
\usepackage{amsmath} 
\usepackage{amssymb}  
\usepackage[style=ieee,url=false,doi=false]{biblatex}
\usepackage[export]{adjustbox}
\usepackage{booktabs}
\usepackage{siunitx}
\usepackage{capt-of}
\usepackage{algorithm}
\usepackage{algpseudocode}
\usepackage{cleveref}
\usepackage{setspace}
\let\labelindent\relax
\usepackage{enumitem}

\usepackage{mathtools}

\newtheorem{defn}{Definition}

\newtheorem{assum}[defn]{Assumption}

\newtheorem{thm}[defn]{Theorem}

\providecommand{\R}{\ensuremath \mathbb{R}}
\providecommand{\Z}{\ensuremath \mathbb{Z}}
\providecommand{\N}{\ensuremath \mathbb{N}}

\newcommand{\ie}{\textit{i.e. }}
\newcommand{\eg}{\textit{e.g. }}

\newcommand{\card}[1]{\left\vert#1\right\vert}

\newcommand{\norm}[1]{\left\Vert#1\right\Vert}

\newcommand{\inv}{^{-1}}

\newcommand{\bigO}{\mathcal{O}}

\newcommand{\defeq}{\vcentcolon=}
\newcommand{\derive}{\mathrm{D}}

\newcommand{\vect}[1]{#1}
\newcommand{\mat}[1]{#1}

\newcommand{\vep}{\varepsilon}
\newcommand{\C}{\mathbb{C}}
\newcommand{\E}{\mathbb{E}} 
\newcommand{\real}{\mathrm{Re}}

\newcommand{\Spec}{\mathrm{Spect}}

\newcommand{\fstate}{\vect{x}}
\newcommand{\rstate}{\vect{p}}
\newcommand{\ctrl}{\vect{u}}
\newcommand{\out}{\vect{y}}

\newcommand{\outrefphase}{z^\star}
\newcommand{\phase}{\vect{\varphi}}
\newcommand{\freq}{\vect{\omega}} 
\newcommand{\ctrltaylor}{\vect{D}}
\newcommand{\ssm}{\vect{W}}
\newcommand{\rdyn}{\vect{R}}
\newcommand{\nnm}{\gamma_{\epsilon}}
\newcommand{\autSSM}{\mathcal{W}(\E)}
\newcommand{\nonautSSM}{\mathcal{W}(\E,\nnm)}

\newcommand{\nfstate}{N}
\newcommand{\nrstate}{n}
\newcommand{\nctrl}{m}
\newcommand{\nout}{o}
\newcommand{\nfreq}{1}
\newcommand{\phasedom}{S^1}
\newcommand{\freqdom}{\R_+}
\newcommand{\nctrltaylor}{\Upsilon}
\newcommand{\nctrlfourier}{\Gamma}

\title{\LARGE \bf
Using Spectral Submanifolds for Nonlinear Periodic Control
}

\author{Florian Mahlknecht$^{1}$, John Irvin Alora$^{1}$, Shobhit Jain$^{2}$,\\Edward Schmerling$^{1}$, Riccardo Bonalli$^{3}$, George Haller$^{2}$, Marco Pavone$^{1}$%
\thanks{
$^{1}$ Department of Aeronautics and Astronautics, Stanford University, Stanford, CA, 94305, USA {\tt\small\{florianm, jjalora, schmrlng, pavone\}@stanford.edu}}
\thanks{
$^{2}$ Institute of Mechanical Systems, ETH Zurich, 8092 Zurich, Switzerland {\tt\small\{shjain, georgehaller\}@ethz.ch}}
\thanks{$^{3}$Laboratory of Signals and Systems, University of Paris-Saclay, CNRS, CentraleSup\'elec, {\tt\small riccardo.bonalli@l2s.centralesupelec.fr}}
}

\begin{document}

\maketitle
\thispagestyle{empty}
\pagestyle{empty}

\begin{abstract}

Very high dimensional nonlinear systems arise in many engineering problems due to semi-discretization of the governing partial differential equations, \eg through finite element methods.
The complexity of these systems present computational challenges for direct application to automatic control.
While model reduction has seen ubiquitous applications in control, the use of nonlinear model reduction methods in this setting remains difficult.
The problem lies in preserving the structure of the nonlinear dynamics in the reduced order model for high-fidelity control.
In this work, we leverage recent advances in Spectral Submanifold (SSM) theory to enable model reduction under well-defined assumptions for the purpose of efficiently synthesizing feedback controllers.

\end{abstract}

\section{INTRODUCTION}

Automatic control of complex, infinite-dimensional systems (\textit{i.e.}, dynamically evolving continua) such as soft robots as well as aircraft and underwater vehicles with coupled fluid-structure interactions remains challenging. 
Model reduction provides a principled approach to reduce model complexity while capturing the essential physics required for controller synthesis.
In optimal control, we are interested in optimizing over a set of control inputs to track a desired trajectory or stabilize around an operating point. 
In these settings, working with the full-order model (FOM) is computationally intractable.
The need to utilize high-fidelity models to control these challenging systems has resulted in significant research efforts towards application of reduced-order models (ROM) for controller design.

In this work we explore recent developments in Spectral Submanifold (SSM) theory \cite{haller2016nonlinear} for nonlinear model reduction and control.
SSMs are the smoothest invariant manifolds that act as nonlinear continuations of the eigenspaces from linearization of a system at a fixed point. 
This nonlinear continuation is tangent to a corresponding spectral subspace of the linearized system.
The additional structure given by this continuation allows us to capture highly-nonlinear behavior outside the vicinity of the linear approximation.
Under certain conditions, these nonlinearities can be approximated arbitrarily well without ever increasing the size of the ROM.

\textbf{Contributions:} Motivated by the established theory on SSMs and their successful application to model reduction \cite{jain2021compute}, we propose the adaptation of SSM theory to automatic control.
In particular, we aim to synthesize optimal, nonlinear, periodic, feedback controllers on the reduced-order SSM to exploit  the computational speed-up while retaining the fidelity of the high-order model.

Our contributions are threefold:

\begin{enumerate}[label=(\roman*)]
    \item We present, for the first time, this novel model reduction technique in the context of control, enabling us to track (quasi-)periodic trajectories. In this setting we restrict ourselves to a special case of trajectory tracking and synthesize optimal control laws that guarantee the existence and persistence of lower-dimensional SSMs on which the true system trajectory lies.
    \item We demonstrate the utility of SSMs in accurately capturing the nonlinear ``slow dynamics" of a system while neglecting its ``fast dynamics". This gives a natural setting in which our control effort is focused on the dynamics that persist.
    \item We illustrate our approach on a pedagogical example, highlight its advantages and disadvantages, and motivate its application to high-dimensional models.
\end{enumerate}


\textbf{Related work: }
Most of the applications of model reduction for control exploit projection-based methods.
They involve a data-driven procedure to identify a linear subspace from simulation rollouts of the FOM.
This approach has been successfully leveraged in literature for real-time control of infinite-dimensional systems using Model Predictive Control (MPC).
The work in \cite{LorenzettiMcClellanEtAl2021} considers the effect of proper orthogonal decomposition (POD) on the closed-loop error dynamics and the authors develop a constraint-tightening scheme to ensure satisfaction of safety constraints in an MPC framework.
While similar works \cite{ghiglieri2014optimal, alla2015asymptotic, altmuller2014model} adopt POD-based MPC schemes for the control of certain PDE-classes, others explore different combinations of subspace identification and optimal control schemes.
In \cite{alla2017error}, the authors investigate the sub-optimality of LQR due to the projection error introduced by POD, while \cite{antil2010domain} considers linear quadratic optimal control using balanced truncation.

Direct application of projection-based methods for nonlinear systems is difficult since evaluation of the nonlinear terms results in a more expensive procedure than evaluating the full model directly, as the change of reference frame involves high-dimensional matrix multiplications \cite{farhat2015structure}.
To overcome this limitation, much of the literature involves construction of locally approximating \textit{linear} ROMs for which standard linear control techniques can be applied.
In \cite{huang2020balanced}, the authors propose an iterative LQR scheme combined with balanced truncation to control the 1D Burger's Equation.
The authors in \cite{TonkensLorenzettiEtAl2021} apply POD in a piecewise-affine fashion by reducing linear approximations of the high-fidelity model.
They then evaluate the nonlinearities through interpolation of the linear approximations and apply an MPC framework to control a soft robot.

While these approaches have been demonstrated to work well on various real-world systems, their performance and theoretical guarantees are limited to linear ROMs \cite{LorenzettiMcClellanEtAl2021}.
For highly nonlinear systems operating in less constrained workspaces, linear ROMs can result in low-fidelity surrogates that exhibit poor closed-loop performance and even instability.
This necessitates the need to capture the structure of the nonlinearities in a more direct way and we propose a new direction for addressing nonlinear model reduction for control.

\textbf{Organization:} In \Cref{sec:preliminaries}
 we introduce our notation and definitions used in this work.
\Cref{sec:problem-statement} defines the optimal control problem where we introduce the tracking error in the periodic orbit of the FOM.
In the methodology in \Cref{sec:methodology}, we describe how we achieve model reduction using SSM in our setting and how we leverage the reduced representation of the dynamics to optimize the tracking error previously introduced in the full-order state-space.
\Cref{sec:validation} showcases the application of the theory on an illustrative example and provides the insights needed for tackling higher-dimensional examples.
We conclude this work in \Cref{sec:conclusions}, highlighting the most promising future avenues.

\section{PRELIMINARIES}
\label{sec:preliminaries}
This section provides the preliminaries that contextualize our approach. We first describe the system dynamics model in Section \ref{section:model} and then define necessary notions in \Cref{sub:spectral-subspace} to lay the groundwork for SSM theory.

\subsection{Notation}
\label{sub:notation}
The set of integers and reals are denoted by $\Z$ and $\R$, with their non-negative counterparts denoted by $\Z_+$ and $\R_+$.
The complex numbers are denoted by $\C$. $\phasedom = \R^\nfreq / \left(2\pi\Z^\nfreq\right)$ represents the circle on the real line.
$C^k$ represents the space of $k$ continuously-differentiable functions and $C^a$ represents the space of analytic functions.
$L^2(V, W)$ is the space of square integrable functions from a complete vector space $V$ to $W$.
$\bigO(\cdot)$ represents the standard big-O notation. $\otimes$ is the tensor product, where $z^{\otimes 3} = z \otimes z \otimes z$.

\subsection{System Model}
\label{section:model}
\subsubsection{Full Order Model}
Consider the following continuous-time, control-affine, nonlinear dynamics with equilibrium point at the origin
\begin{align}
\label{eq:FOM}
\begin{cases}
\begin{aligned}
    \dot{\fstate}(t) &= \mat{A} \fstate(t) + \vect{f}_0(\fstate(t)) + \sum_{i=1}^m \vect{f}_i(\fstate(t))\, \ctrl_i(t), \\
    \out(t) &= \mat{H} \fstate(t),
\end{aligned}
\end{cases}
\end{align}
where the state $\fstate(t) \in \R^\nfstate$ is high-dimensional, \ie $\nfstate $ is large; $A \in \R^{\nfstate\times\nfstate}$ is the stability matrix; $f_0(x)$ are the nonlinearities of the uncontrolled system; $f_i : \R^\nfstate \rightarrow \R^\nfstate$ for $i=1,...,m$ are nonlinear functions that describe the state-dependence of the control effort via an $\nctrl$-dimensional control input $\ctrl(t) \in \R^\nctrl$; the observed output of the system is denoted as $\out(t) \in \R^\nout$; and $H \in \R^{\nout \times \nfstate}$ is the selection matrix of output variables, where $\nout \ll \nfstate$. 
In this work, we assume that the performance and output variables are the same and that they are perfectly observable. 
We introduce the following assumption on the form of $A$.

\begin{assum}
\label{assum:Astability}
$A$ is negative definite, \ie $A \prec 0$.
\end{assum}
In other words, we assume that the origin $\bar{\fstate} = \vect{0}$ is a locally asymptotically stable equilibrium point.
Many physical systems and phenomena of interest such as soft robots and fluid structure interactions satisfy this assumption (possibly up to a shift in origin).

In addition, we introduce the following assumption on the form of $f_i$.

\begin{assum}
\label{assum:analytic}
The functions $f_0,\dots,f_m \in C^a$.
\footnote{We make this assumption for ease of exposition. In general, the right-hand side of System~\eqref{eq:FOM} can have finite smoothness, infinite smoothness, or be analytic; correspondingly the spectral submanifold defined in Section~\ref{sec:SSMprelims} is as smooth as the right-hand side.}
\end{assum}
We remark that this assumption is not particularly limiting since many physical systems (\eg soft robots) generically satisfy this assumption and we are only interested in controlling smooth behavior.

\subsection{Spectral Subspace}
\label{sub:spectral-subspace}

Consider the uncontrolled part of System~\eqref{eq:FOM}
\begin{align}
\label{eq:FOM_uncontrolled}
    \dot{\fstate}(t) &= \mat{A} \fstate(t) + \vect{f}_0(\fstate(t)),
\end{align}
whose linearization around the origin is given by
\begin{align}
\label{eq:FOM_uncontrolled_lin}
    \dot{\fstate}(t) &= \mat{A} \fstate(t).
\end{align}

For any eigenvalue $\lambda_j$ of $A$, there exists an eigenspace $\E_j \subset \R^\nfstate$ spanned by the (generalized) eigenvectors of $A$.
These eigenspaces are invariant subspaces of the linearized system~\eqref{eq:FOM_uncontrolled_lin}.

\defn{A spectral subspace $\E_{j_1, ..., j_\nrstate}$ of System~\eqref{eq:FOM} is defined as the direct sum of an arbitrary collection of eigenspaces of $A$ \ie}
\begin{align*}
    \E \defeq \E_{j_1, ..., j_\nrstate} = \E_{j_1} \oplus \E_{j_2} \oplus ... \oplus \E_{j_\nrstate}.
\end{align*}
By linearity of System~\eqref{eq:FOM_uncontrolled_lin}, any spectral subspace is an invariant subspace of $A$. 
In projection-based methods, ROMs are constructed by projecting the dynamics onto a nested hierarchy of the slowest $k$ spectral subspaces \ie $\E^k = \E_{1,...,k}$ where $\E^1 \subset \E^2 \subset \E^3 \subset \dots \subset \E^k$ and $k \ll N$. 

However, such projections of the governing equations to spectral subspaces can be guaranteed to work only for linear systems and do not capture the effects of the nonlinear terms and control inputs of the FOM.
To find a faithful reduction of System~\eqref{eq:FOM}, it is necessary to reason about how the additional nonlinear terms and time dependent forcing influence the structure of the spectral subspace. To this end, we propose using SSMs and their reduced dynamics for reducing the following nonlinear control problem.

\section{PROBLEM STATEMENT}
\label{sec:problem-statement}

In this section we provide a formal problem definition of the full-order, periodic optimal control problem in Section \ref{ROMPOCP}.

\subsection{Periodic Optimal Control Problem}
\label{ROMPOCP}

In this work we design periodic orbits, minimizing the mean distance to some desired trajectory $\outrefphase(\freq\,t) \in \R^o$, where $\freq \in \freqdom$ is the frequency of the reference trajectory.
Our approach is to formulate the following optimal control problem
\begin{align}
\label{eq:PeriodicOCP}
    \min_{u(\cdot)} \quad & \frac{1}{T} \int_{0}^{T} \norm{\outrefphase(\freq\, t) - \out(t)}_2 dt \nonumber \\
    \text{subj. to}& \quad 
    \text{System~\eqref{eq:FOM}} \\
    & \quad \fstate(0) = \fstate(T) \nonumber .
\end{align}

In Equation~\eqref{eq:PeriodicOCP}, we minimize over a class of periodic feedback control laws of the form
\begin{align}
    \ctrl(t) = \kappa(\out(t), \freq t),
\end{align}
where $\kappa \in L^2(\R^\nout \times \phasedom, \R^\nctrl)$.
For ease of notation, throughout the rest of the paper we denote $\phase = \freq t$.

Informally, we minimize the mean-squared trajectory error between our system's periodic orbit and the desired trajectory, after its fast dynamics have sufficiently decayed.
We emphasize that in this work, we are interested in synthesizing control laws that neglect transients and control for a periodic orbit.

We remark that while we consider the case of periodic control laws in this paper, our approach generalizes to the quasi-periodic setting. 

\section{METHODOLOGY}
\label{sec:methodology}

\subsection{Spectral Submanifold Preliminaries}
\label{sec:SSMprelims}
An SSM serves as the unique nonlinear continuation of a nonresonant spectral subspace $\E$ for the nonlinear system~\eqref{eq:FOM_uncontrolled} and is defined as follows~\cite{haller2016nonlinear}.

\defn{An autonomous SSM $\autSSM$, corresponding to a spectral subspace $\E$ of the operator $A$ is an invariant manifold of the nonlinear system~\eqref{eq:FOM_uncontrolled} such that}
\begin{enumerate}
    \item $\autSSM$ is tangent to $\E$ at the origin and has the same dimension as $\E$,
    \item $\autSSM$ is strictly smoother than any other invariant manifold satisfying condition 1 above. 
\end{enumerate}
A slow SSM is associated to a spectral subspace containing the slowest decaying eigenvectors of the linearized system. Slow SSMs are ideal candidates for model reduction as typical nearby full system trajectories are exponentially attracted towards these manifolds and synchronize with the slow dynamics on such SSMs.


We synthesize such a controller by focusing on controlling the reduced dynamics along a slow SSM.
As the full system trajectories quickly and automatically synchronize with the dynamics on the slow SSM, we envision a minimal control effort arising from our synthesized controller on the slow SSM.
Hence, we assume a small control input by rescaling the control terms in system~\ref{eq:FOM} by a small scalar parameter $\vep>0$ as

\begin{align}
\label{eq:FOM_epsilon}
\begin{cases}
\begin{aligned}
    \dot{\fstate}(t) &= \mat{A} \fstate(t) + \vect{f}_0(\fstate(t)) + \epsilon g(x(t),\freq t), \\
    g(x(t),\freq t) &= \sum_i^m \vect{f}_i(\fstate(t))\, \kappa_i(\mat{H} \fstate(t), \freq t),
\end{aligned}
\end{cases}
\end{align}
where the control input $\kappa_i(\out(t), \freq t)$ has periodic time-dependence with frequency $\freq$ for all $i=1,\dots,m$. 

In this non-autonomous setting of periodic control, SSMs are envisioned similarly to the autonomous setting and the role of the fixed point is taken over by the periodic orbit $\nnm$ created by the small-amplitude control force. A nonautonomous, time-periodic SSM $\nonautSSM$ is then a fibre bundle that perturbs smoothly from the vector bundle $\nnm\times \E$ under the addition of the nonlinear and control terms in System~\eqref{eq:FOM_epsilon}. Hence, $\nonautSSM$ is $\frac{2\pi}{\freq}$-periodic in time.

\defn{A time-periodic SSM $\nonautSSM$, corresponding to a spectral subspace $\E$ of the operator $A$ is an invariant manifold of the nonlinear system~\eqref{eq:FOM_epsilon} such that}
\begin{enumerate}
    \item $\nonautSSM$ is a subbundle of the normal bundle $N\nnm$ of the periodic orbit $\nnm$, satisfying $\dim \nonautSSM = \dim \E + 1$,
    \item $\mathcal{W}(\E)$ perturbs smoothly from the spectral subspace $\E$ of the linearized system  under the addition of nonlinear and control terms in System~\ref{eq:FOM_epsilon}.
    \item $\nonautSSM$ has strictly more continuous derivatives along $\nnm$ than any other invariant manifold satisfying conditions 1 and 2 above. 
\end{enumerate}

For any spectral subspace $\E$, the \emph{absolute spectral quotient}~\cite{haller2016nonlinear} is defined as
\begin{align}
\label{eq:spec_qt}
    \Sigma(\E) &= \mathrm{Int} \left[  \frac{\min_{\lambda\in \Spec(A)}\real \lambda}{\max_{\lambda\in \Spec(A|_\E)}\real \lambda}  \right].
\end{align}
This spectral quotient measures the fastest decay exponent outside the spectral subspace $\E$ relative to the slowest decay exponent within $\E$. It is crucial for determining the smoothness class of invariant manifolds in which the SSM uniquely exists. A high-value of the spectral quotient indicates a high-degree of overlap between invariant manifolds tangent to $\E$ at the origin, which is desirable for model reduction over slow SSMs.

For a small-enough control effort, the following theorem guarantees the existence of a time-periodic SSM, whose reduced dynamics provides us an exact nonlinear reduced-order model for control synthesis.

\thm{\label{thm:existence}Consider a spectral subspace $\E$ with $\dim \E = \nrstate$ and its associated eigenvalues (counting multiplicities) listed as $\lambda_1,\dots,\lambda_{\nrstate}$. Assume that the low-order nonresonance conditions 
\begin{equation}
    \label{eq:nonres}
    \sum_{j=1}^{\nrstate} m_j \real \lambda_j \ne \real \lambda_l, \quad \lambda_l\not\in \Spec(A|_{\E}),\quad 2\le \sum_{j=1}^{\nrstate} m_j\le \Sigma(\E),
\end{equation}
hold for all eigenvalues $\lambda_l$ of $A$ that lie outside the spectrum of $A|_{\E}$ with  $m_j \in \N$ and that Assumptions~\eqref{assum:analytic} and \eqref{assum:Astability} are satisfied.

Then the following holds:
\begin{enumerate}
    \item There exists a time-periodic SSM, $\nonautSSM$ for system~\eqref{eq:FOM_epsilon} that depends smoothly on the parameter $\epsilon$ and is unique in the class of $C^{\Sigma(\E) +1}$ invariant manifolds.
    \item $\nonautSSM$ can be viewed as an embedding of an open set $\mathcal{U}$ into the state space of System~\eqref{eq:FOM_epsilon} via the map
    \begin{equation}
        \ssm_\vep(\rstate,\phase):\mathcal{U}\subset \C^{\nrstate}\times \phasedom \to \R^{\nfstate},
    \end{equation}
    with the periodic phase variable $\phase\in \phasedom$.
    \item There exists a polynomial function with respect to $\rstate$, $\rdyn_\vep(\rstate,\phase):\mathcal{U}\to \C^{\nrstate}$ satisfying the invariance equation
    \begin{align}
        A\ssm_\vep(\rstate,\phase) + f_0(\ssm_\vep(\rstate,\phase)) + \epsilon g(\ssm_\vep(\rstate,\phase),\phase) = \nonumber \\\derive_p\ssm_\vep(\rstate,\phase)\rdyn_\vep(\rstate,\phase) + \derive_{\phase}\ssm_\vep(\rstate,\phase)\freq,
        \label{eq:invariance}
    \end{align}
    such that the reduced dynamics on the SSM is given by 
    \begin{equation}
    \label{eq:ROM}
        \dot{\rstate} = \rdyn_\vep(\rstate,\phase).
    \end{equation}
\end{enumerate}
}

\begin{proof}
This is a restatement of Theorem 4 in~\cite{haller2016nonlinear} in our setting, which is deduced from the abstract results on whiskers of invariant tori in \cite{Haro2006}.
\end{proof}

\subsection{Model reduction using SSM}
Theorem~\ref{thm:existence} allows us to approximate $\ssm_\vep(\rstate,\phase),\rdyn_\vep(\rstate,\phase)$ in a neighbourhood of the origin as a Taylor expansion in the parametrization coordinates $\rstate$ with coefficients that depend periodically on the phase variable $\phase$. These periodic cofficients can be further Fourier-expanded resulting in Taylor-Fourier series for $\ssm_\vep(\rstate,\phase),\rdyn_\vep(\rstate,\phase)$. This means that the SSM and its reduced dynamics can be approximated arbitrarily well without ever increasing the dimension of $\E$. This is a highly desirable property for control since it enables one to faithfully capture the essential nonlinearities in the dynamics without increasing the dimensionality of the model.

As detailed in \cite{jain2021compute}, the solution of the invariance Equation~\eqref{eq:invariance} can be efficiently accomplished by solving the mappings $\ssm_\vep$ and $\rdyn_\vep$ with the ansatz

\begin{align}
    \ssm_\vep(\rstate, \phase) &= \ssm_0(\rstate) + \vep \ssm_1(\rstate, \phase) + \bigO(\vep^2), \\
    \rdyn_\vep(\rstate, \phase) &= \rdyn_0(\rstate) + \vep \rdyn_1(\rstate, \phase) + \bigO(\vep^2),
\end{align}
where the autonomous terms with $\vep = 0$ are expressed as multivariate Taylor expansions:

\begin{align}
    \ssm_0(\rstate) &= \sum_{j \ge 0} \ssm_{0,j} \rstate^{\otimes j}, \\
    \rdyn_0(\rstate) &= \sum_{j \ge 0} \rdyn_{0,j} \rstate^{\otimes j},
\end{align}
with the unknown coefficients $\ssm_{0,j}$, $\rdyn_{0,j}$ being $(j+1)$-tensors.
The $\mathcal{O}(\vep)$ terms are expanded via a Taylor-Fourier series as
\begin{align}
\label{eq:polynomial-eps-ssm}
    \ssm_1(\rstate,\phase) &= \sum_{j \ge 0} \ssm_{1,j}(\phase) \rstate^{\otimes j}, \quad \ssm_{1,j}(\phase) = \sum_{\vect{h} \in \Z^\nfreq}\ssm_{1,j,h} e^{i \langle \vect{h}, \phase \rangle} \\
\label{eq:polynomial-eps-rdyn}
    \rdyn_1(\rstate,\phase) &= \sum_{j \ge 0} \rdyn_{1,j}(\phase) \rstate^{\otimes j}, \quad \rdyn_{1,j}(\phase) = \sum_{\vect{h} \in \Z^\nfreq}\rdyn_{1,j,h} e^{i \langle \vect{h}, \phase \rangle},
\end{align}
with $\ssm_{1,j,h}, \rdyn_{1,j,h}$ denoting unknown Taylor-Fourier coefficients at degree $j$ and harmonic $h \in \N$. As detailed in~\cite{jain2021compute}, these unknown coefficients are determined by solving the invariance equation~\eqref{eq:invariance} in a recursive manner, where each recursion involves the solution of a linear system. These computations have been automated and demonstrated on nonlinear finite-element based applications featuring more than 100,000 degrees of freedom~\cite{jain2021compute}. SSMTool, an open-source implementation of this procedure is available at~\cite{SSMTool}.

\subsection{Exploiting the ROM for offline optimization}

We consider a generic periodic feedback control law, expressible through a truncated Taylor-Fourier series:

\begin{align}
\label{eq:ctrl-ansatz-taylor}
    \kappa(\out, \phase)  =  \sum_{j = 0}^{\nctrltaylor} \ctrltaylor_j(\phase) \out^{\otimes j},
\end{align}
where $\nctrltaylor \in \N$ is the finite truncation order of the Taylor series. $\ctrltaylor_j(\phase)$ is a tensor of order $j+1$ and dimension $\nctrl$, \ie $\ctrltaylor_0(\phase) \in \R^\nctrl$.
The coefficients are individually determined by the following (truncated) Fourier series:

\begin{equation}
\label{eq:ctrl-ansatz-fourier}
    \ctrltaylor_j(\phase) = \sum_{\vect{h} \in \mathbb{H} \subset \Z^\nfreq} \ctrltaylor_{j,\vect{h}} e^{i \langle \vect{h} , \phase \rangle},
\end{equation}
where $\nctrlfourier = \card{\mathbb{H}}$ is the finite truncation order of the Fourier series.
This allows us to consider the controller family $\kappa_{\ctrltaylor_{j,\vect{h}}}$ generated by all possible realizations of the parameters $\ctrltaylor_{j,\vect{h}} \in \R^{\nctrlfourier\times\nctrl\times\nout^j}, \, j \in \{0,\dots,\nctrltaylor\}$.
The number of parameters to optimize is therefore $n_p = \sum_{j = 0}^\nctrltaylor\, \nctrlfourier\, \nctrl\, \nout^j$.

As a consequence, the previously derived mappings of the $\bigO(\vep)$-perturbed SSM and its reduced dynamics in \cref{eq:polynomial-eps-ssm,eq:polynomial-eps-rdyn} are now dependent on the control parameters, \ie $\ssm_\vep(\rstate, \phase, \ctrltaylor_{j,\vect{h}})$, $\rdyn_\vep(\rstate, \phase, \ctrltaylor_{j,\vect{h}})$.

We exploit this reduced order presentation, to find the optimal parameters $\ctrltaylor_{j,\vect{h}}^\star$ in an offline optimization procedure.
The ROM optimization formulation reads:

\begin{align}
\label{eq:reduced-optimization}
    \min_{\out(\cdot), \> \ctrltaylor_{j,\vect{h}}} \quad \frac{1}{T} & \int_{0}^{T} \norm{\outrefphase(\freq\, t) - \out(t)}dt \\
    \text{subj.} \quad \dot{\rstate} &= \rdyn_\vep(\rstate, \phase, \ctrltaylor_{j,\vect{h}}) \nonumber \\
    \out &= \mat{H} \, \ssm_\vep(\rstate, \phase, \ctrltaylor_{j,\vect{h}}) \nonumber \\
    \out(0) &= \out(T) \nonumber .
\end{align}
We remark that $\rstate \in \R^\nrstate$, $\nrstate \ll \nfstate$. 
Hence, this optimization problem is much more tractable than Problem \ref{eq:PeriodicOCP}, motivating the construction of the reduced model.

\subsection{Summary}
\label{sec:summary}

We summarize our method in Algorithm~\ref{alg:nlmr-periodic-ctrl-ssm}.
As an input we process the system matrices of System~\eqref{eq:FOM} with the asymptotically stable fixed point shifted to the origin.
$f_i(\fstate)$ are defined by supplying the coefficients of their respective multivariate Taylor expansions.

\begin{algorithm}[h]
\caption{Periodic control with SSM}\label{alg:nlmr-periodic-ctrl-ssm}
\begin{algorithmic}[1]
\Require
\hspace*{-1em}\begin{itemize}
    \item System~\eqref{eq:FOM}: $\mat{A} \in \R^{\nfstate \times \nfstate}, \, \vect{f}_i: \R^\nfstate \to \R^\nfstate,\, \mat{H} \in \R^{\nout \times \nfstate}$
    \item $\outrefphase(\phase)$
\end{itemize}
\setstretch{1.5}
\Ensure $ \mat{A} \fstate(0) + \vect{f}_0(\fstate(0)) + \sum_i^m \vect{f}_i(0)\, \ctrl(0) = \vect{0}$
\State $\lambda_i, \vect{v}_i \gets$ \Call{SpectralDecomposition}{$\mat{A}$}
\State $E^n$ $\gets \bigoplus_{k \in \{j_1 \dots j_\nrstate\}} v_k $
\Comment{Pick the $\nrstate$ slowest dynamics}
\State Define form of $\kappa_{\ctrltaylor_{j,\vect{h}}}(\out, \phase)$ \Comment{dependent on $\ctrltaylor_{j,\vect{h}}$}
\State $\ssm_{\vep,\ctrltaylor_{j,\vect{h}}}, \rdyn_{\vep,\ctrltaylor_{j,\vect{h}}} \gets$ \Call{ComputeSSMCoeffs}{$\mat{A}, f_i, E^n, \ctrltaylor_{j,\vect{h}}$}
\State $\ctrltaylor_{j,\vect{h}}^\star \gets$ \Call{Optimize} {\cref{eq:reduced-optimization}}
\State Apply feedback law $u(t) = \kappa_{\ctrltaylor_{j,\vect{h}}^\star}(\out(t), \phase)$ to System \ref{eq:FOM}
\end{algorithmic}
\end{algorithm}

The control law is defined by picking the expansion order $\nctrltaylor$ and the $\nctrlfourier$ integer combinations of desired frequency components, which determines the number of parameters $n_p$ that we optimize over.

Once optimal parameters for following the trajectory $\outrefphase(\phase)$ are found, we apply our optimal feedback periodic control law to the FOM in System~\ref{eq:FOM}.
Assuming the designed periodic orbits are stable, this control strategy guarantees that our system trajectories will asymptotically converge to the $\varepsilon$-perturbed SSM containing this orbit.

\section{VALIDATION}
\label{sec:validation}

\subsection{Overdamped Pendulum Dynamics}
To illustrate the principles of SSM theory for control, we consider the simple example of an overdamped pendulum providing a two-dimensional spectral subspace associated with two distinct stable eigenvalues. 
We denote the state space variables with $x_1 = \theta,\, x_2 = \dot{\theta}$.
The dynamics of the system are then given as
\begin{equation}
\label{eq:pendulum_dyn}
  \begin{cases}
    \dot{x}_1(t) = x_2(t) , \\
    \dot{x}_2(t) = -\frac{b}{m\ell^2}\,x_2(t) - \frac{g}{\ell} \sin x_1(t) + \varepsilon \frac{1}{m\ell^2} u(t).
  \end{cases}
\end{equation}

Considering the fixed point to be at the origin (corresponding to the pendulum in the downward position with no motion), we convert the system to the form denoted in Equation~\eqref{eq:FOM} by splitting it into a linear part $A$, the nonlinear part $f_0(x)$, and control-affine part:
\begin{equation}
\dot{x}(t) = 
\underbrace{\begin{bmatrix}
0 & 1 \\
-\frac{g}{\ell} & -\frac{b}{m\ell^2}
\end{bmatrix}}_{{A}} x(t) + \underbrace{\begin{bmatrix}
0 \\
\frac{g}{\ell} x_1(t) - \frac{g}{\ell}\sin x_1(t)
\end{bmatrix}}_{{f_0}({x(t)})}+ \varepsilon \begin{bmatrix}
0 \\
\frac{1}{m\ell^2}
\end{bmatrix} u(t) .
\label{eq:pendulum-1st-order}
\end{equation}
Note that $\epsilon$ makes explicit that the magnitude of $u(t)$ should be \textit{moderate}; we provide more insight on this later in this section. 
In our experiments, we set $\epsilon = 1$.

\begin{figure}
\parbox[t]{.4\linewidth}{\null
  \centering
  \includegraphics[height=1in,trim=11 15 10 0,clip]{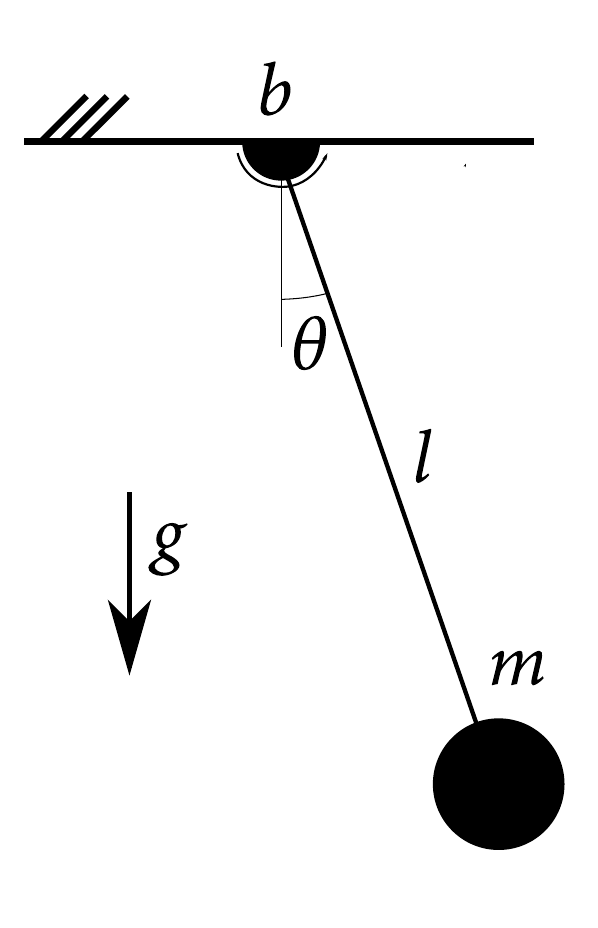}%
  \captionof{figure}{Pendulum illustration\label{fig:pendulumfig}}%
}
\parbox[t]{.6\linewidth}{\null
\centering
  \vskip-\abovecaptionskip
  \captionof{table}[t]{Pendulum parameters\label{tab:pendulumparams}}%
  \vskip\abovecaptionskip
  \resizebox{!}{.45in}{%
    \begin{tabular}{crl}
  \toprule
  \textbf{Parameter} & \textbf{Value} & \textbf{Unit} \\
  \midrule
     m & 1 & \si{kg}\\
     $\ell$ & 1 & \si{m}\\
     g & 9.81 & \si{m/s^2}\\
     b & 35 & \si{Nm s/rad}\\
  \bottomrule
  \end{tabular}}
}
\end{figure}

\begin{figure}
    \centering
    \includegraphics[width=.8\linewidth]{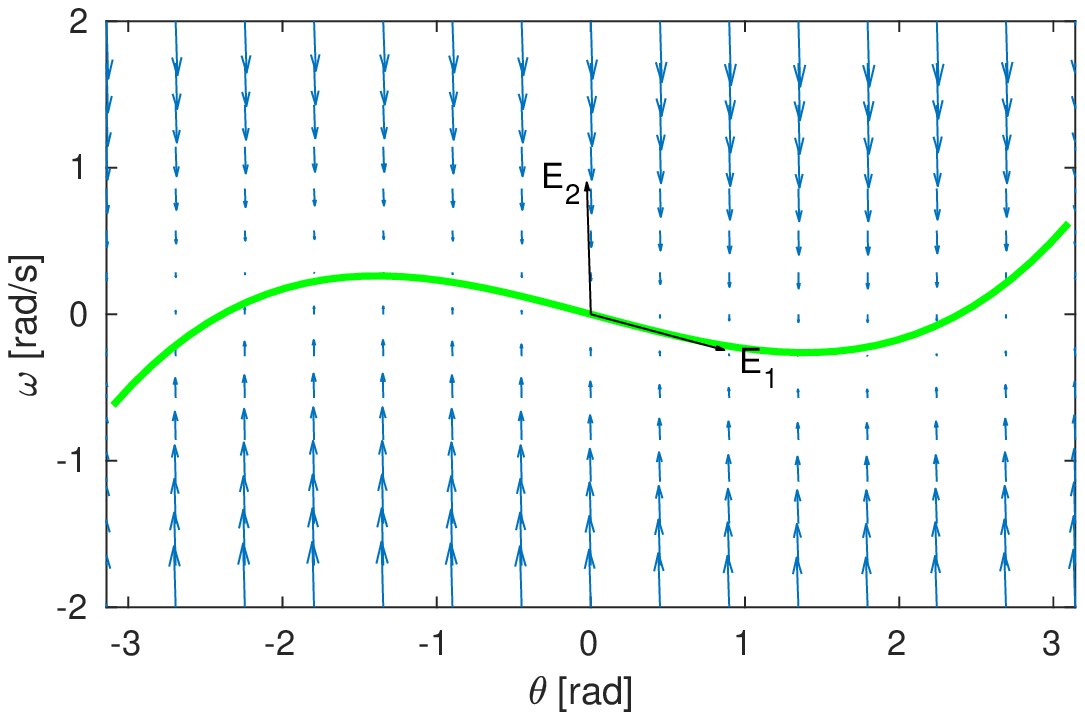}
    \caption{Overdamped pendulum phase portrait with SSM. $\E_1$ represents the slow spectral subspace to which we attach our SSM (depicted in green) while $\E_2$ represents the fast spectral subspace commensurate with the dynamics that converge quickly to the SSM.}
    \label{fig:pendulum-ssm}
\end{figure}

The stability matrix $A$ has two eigenvalues $0 > \lambda_1 > \lambda_2$ with corresponding eigenvectors $\vect{v}_1$ and $\vect{v}_2$.
We pick the spectral subspace spanned by $\vect{v}_1$ (\ie $\E_1 = \text{lin}(\mathbf{v}_1)$) which corresponds to the slowest converging mode in $A$.
The slow and fast spectral subspaces $\E_1$ and $\E_2$, respectively, and the attached SSM to $\E_1$ are shown in Figure~\ref{fig:pendulum-ssm}.
Lastly, using Equation~\eqref{eq:spec_qt}, we compute the spectral quotient to be $\sigma(\E_1) = \mathrm{Int} \left[ \frac{\lambda_2}{\lambda_1}\right] = 122$. Recalling Theorem~\eqref{thm:existence}, we verify that the non-resonance conditions are met. 

We consider periodic state feedback controllers for the pendulum of the following form, defined in terms of coefficients $\vect{u_p} \defeq \begin{bmatrix} u_{p1} \dots  u_{p6}\end{bmatrix}$

\begin{align}
\label{eq:periodicinput}
    \kappa_{\vect{u_p}}(\vect{x},\,\phase) \defeq &u_{p1} + u_{p2} \cos\phase + u_{p3} \sin\phase  \\
    &+ x_1 \left( u_{p4} + u_{p5} \cos\phase + u_{p6} \sin \phase \right) \nonumber
\end{align}

\subsection{SSM Derivation}
To compute an analytic expression of the SSM, we use graph-style parametrization.
Therefore, we first perform a change of basis into the spectral coordinates, \ie
\begin{align*}
  \vect{x} = \mat{T} \vect{\xi} = \begin{bmatrix}\vect{v}_1 & \vect{v}_2\end{bmatrix} \begin{bmatrix}\xi_1 \\ \xi_2\end{bmatrix}.
\end{align*}

We express the SSM as a function over $\E_1$, \ie $\xi_2 = h(\xi_{1}, \phase)$.
Hence, the mapping back to the FOM state space, is given by:
\begin{equation}
\vect{x} =  \ssm(\xi_1) = \mat{T} \begin{bmatrix}
  \xi_1 \\
  h(\xi_1, \phase)
\end{bmatrix}
\end{equation}

Using the transformation $\vect{\dot{\xi}} = \mat{T}\inv\, \vect{f}(\mat{T}\vect{\xi})$ and taylor expanding the non-polynomial, nonlinear terms (denoted $\tilde{f}_{nl}$) around the origin, we obtain the following dynamics
\begin{align*}
\dot{\xi} =
\Lambda \xi + \tilde{f}_{nl}(\xi_1, \xi_2)
+ \vep \mat{T}\inv \begin{bmatrix}
0 \\
\frac{1}{m\ell^2}
\end{bmatrix} \kappa_{\vect{u_p}}(\phase,\ssm(\xi_1)).
\end{align*}

Denoting the dynamics for
$\dot{\xi}_1$ as
$g_1(\xi, \phase)$ and
$\dot{\xi}_2$ as
$g_2(\xi, \phase)$,
we state the invariance equation as
\begin{align}
  \label{eq:pendulum-invariance}\left. g_2 \right|_{\xi_2 = h(\xi_1,\phase)} =
  \left. \derive_{\xi_1} h(\xi_{1}, \phase)\, g_1\right|_{\xi_2 = h(\xi_1, \phase)} + \derive_{\phase} h(\xi_{1},\phase)\, \freq .
\end{align}

The right hand side is given by the derivative in time of our SSM parametrization $h(\xi_{1}, \phase)$, similarly to Equation~(\ref{eq:invariance}).

We solve this invariance equation with the ansatz:
\begin{align*}
  h(\xi_{1}) &= c_1\,\xi_1^2 + c_2\,\xi_1^3
  + \epsilon h_1(\xi_1, \phase) , \\
  h_1(\xi_1,\phase) &= c_3 + c_4 \cos\phase + c_5 \sin\phase + c_6\, \xi_1 \cos\phase + c_7\, \xi_1 \sin\phase,
\end{align*}

By coefficient comparison we determine $\begin{bmatrix}c_1 \dots c_7\end{bmatrix}$ as a function of $u_p$.
In this way we obtain $h_{u_p}(\xi_{1}, \phase)$, representing the perturbed SSM due to the parametric forcing by the periodic feedback controller.

Let us define $p \defeq \xi_1$.
Then, the reduced dynamics of the full system in Equation~\eqref{eq:pendulum_dyn} is represented on the SSM as
\begin{align}
\label{eq:reduced-pendulum}
    \dot{p} = g_1(p, h_{u_p}(p, \phase))
\end{align}

\begin{figure}
    \centering
    \includegraphics[width=\linewidth]{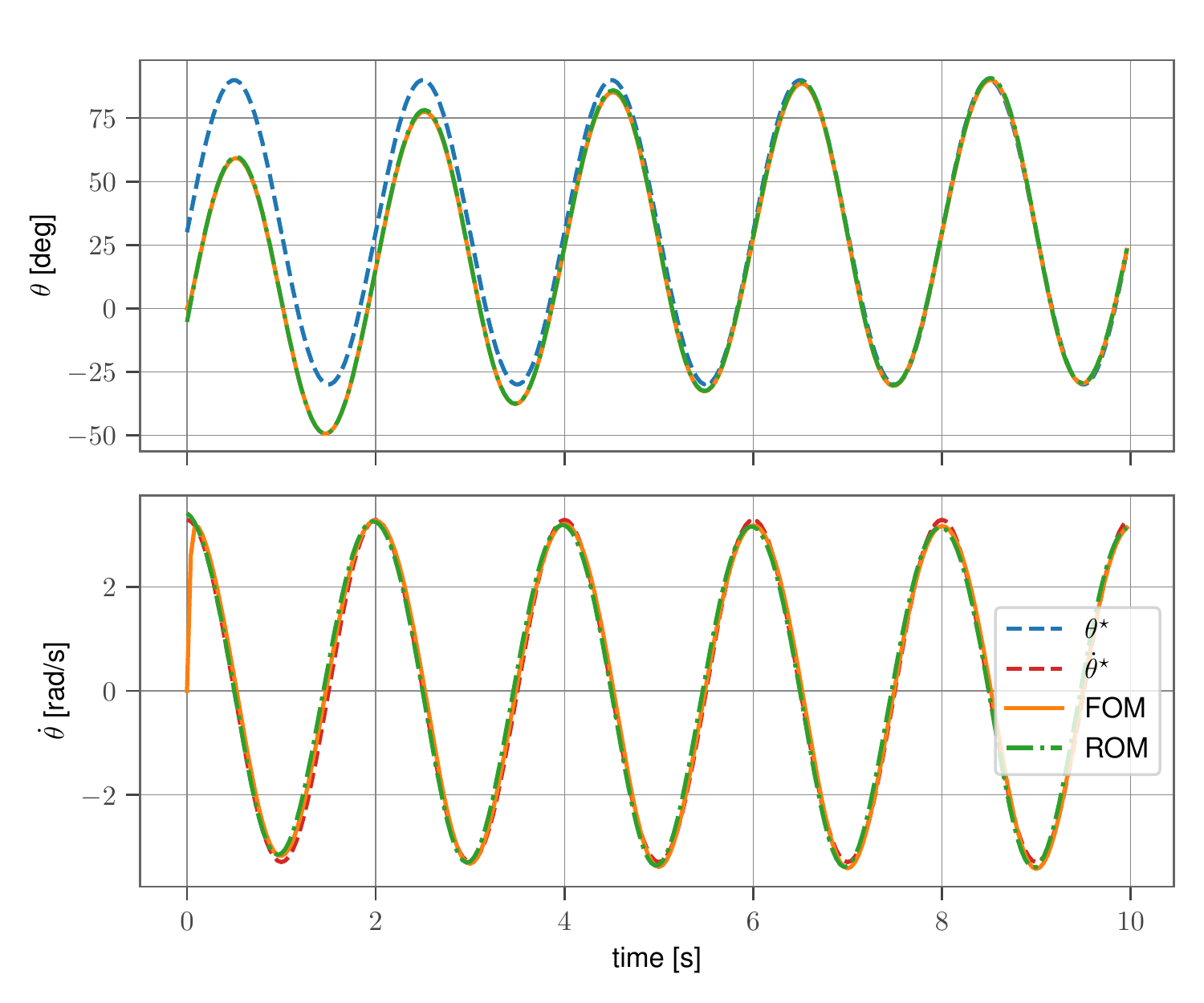}
    \caption{FOM and ROM tracking performance. The top figure represents tracking performance for the desired $\theta^\star(t)$. The bottom figure represents the untracked state $\dot{\theta}$. Simulation is carried out for five time-periods \ie $t_f = 5T$. The steady-state error between the closed-loop system and the desired trajectory is RMSE($\theta$) = $0.49$ deg and RMSE($\dot{\theta}$) = $0.19$ rad/s.}
    \label{fig:theta-tracking-error}
\end{figure}

\begin{figure}
    \centering
    \includegraphics[width=\linewidth]{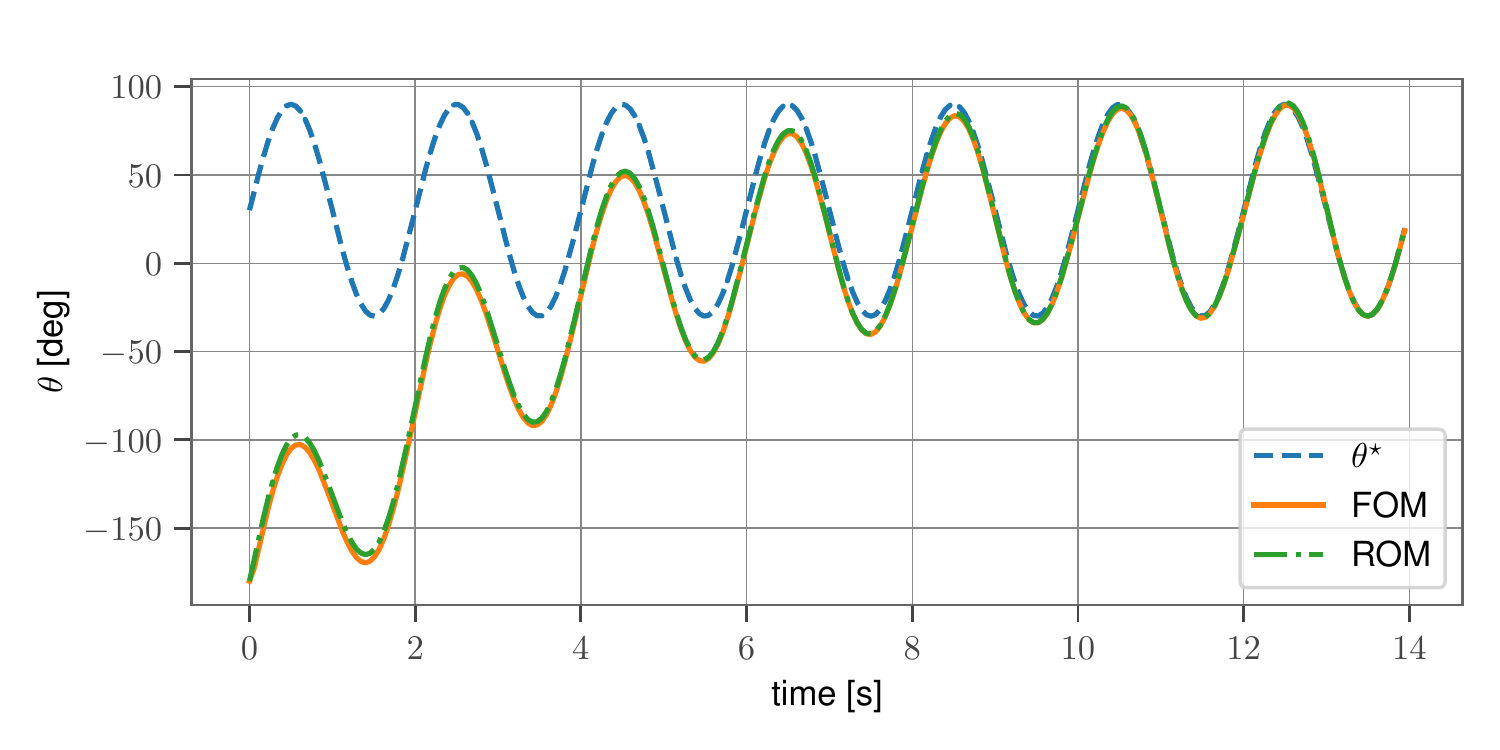}
    \caption{Convergence to the periodic orbit from  $\theta_0$ = -180 deg. The ROM on the SSM perfectly captures the nonlinearities even when far away from the linearization point, \ie the origin.}
    \label{fig:pendulum-far-x0}
\end{figure}

\subsection{Controller Performance}

In this section we hope to showcase the predictive capability of SSMs for control synthesis as well as give insights on the limitations of our approach. 
We motivate the extension of this work to higher-dimensional problems which satisfy the conditions set forth in our pedagogical experiment.

As shown in Algorithm~\ref{alg:nlmr-periodic-ctrl-ssm}, we first compute the reduced-order representation in Equation~\eqref{eq:reduced-pendulum}.
Our desire is to track the trajectory $\theta^\star(\freq t) = 30 + 60 \sin(\freq t)$ (in degrees), where $\freq = \pi$ (rad/s).
We simulate the trajectory and dynamics for five time periods, ensuring that the dynamics in Equation~\eqref{eq:pendulum_dyn} achieves its periodic orbit.
We then compute the optimal parameters of our periodic feedback control law $u(t)$ in Equation~\eqref{eq:periodicinput} by solving Problem~\ref{eq:reduced-optimization}.
We use the CMA-ES optimization algorithm \cite{hansen2006cma} implemented in the KORALI framework \cite{martin2021korali}.
After two time-periods, System~\eqref{eq:pendulum_dyn} achieves its periodic orbit, hence we set $t_1 = 2T$.

Figure~\ref{fig:theta-tracking-error} shows the closed-loop tracking performance for both the reduced model trajectory and the full order trajectory with the optimal coefficients from the ROM optimization.
In other words, we compare the full system's evolution $\fstate(t)$ with the reduced system's evolution $\ssm(\rstate(t))$.
As expected, after the small initial transient, the full system trajectory is quickly attracted to the periodic orbit induced by our controller.
There is no noticeable difference between FOM and ROM trajectories, meaning that the FOM trajectory lies as expected on the SSM and its $\varepsilon$-perturbation captures the periodic motion well.

In Figure~\ref{fig:pendulum-far-x0} we perform the same experiment, but we initialize the pendulum at $\theta=-180$ deg.
Despite the large distance to the origin and the significant nonlinearities, the ROM still evolves in the same manner as the full-order model.
Furthermore, the system converges to the periodic reference trajectory, as desired.

\begin{figure}
    \centering
    \includegraphics[width=.8\linewidth]{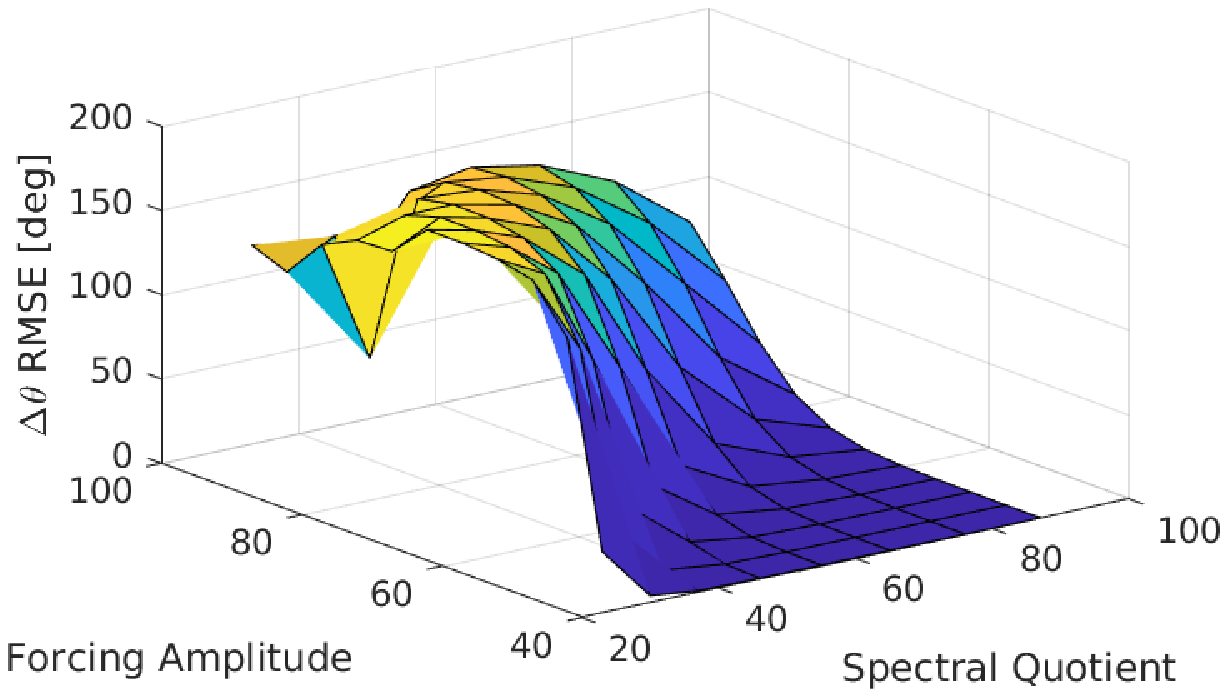}
    \caption{Error comparison between FOM and ROM for different spectral quotients and forcing amplitudes. As expected, a higher forcing amplitude leads to destruction of the manifold, since the $\vep$-order perturbation is not small anymore. The allowed scale of $\bigO (\vep)$ is driven by the spectral quotient, and in particular for our pendulum example, the damping coefficient.}
    \label{fig:ssm-robustness}
\end{figure}
\begin{figure}
    \includegraphics[keepaspectratio, width=8.5cm,clip,trim=0 0 0 0 ]{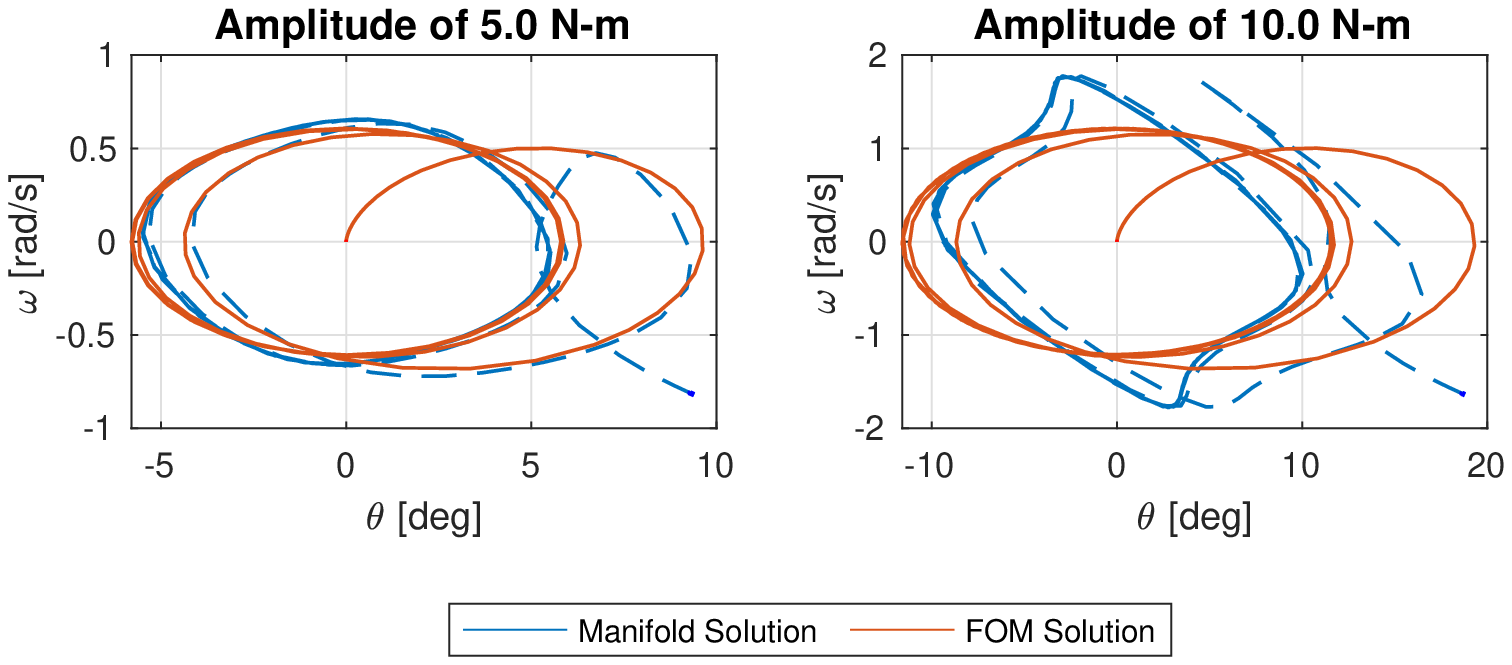}
    \caption{Periodic orbits of a pendulum ($b = 7$) with relatively low spectral quotient ($\sigma(E) = 2$) versus increasing torque amplitude. This explicitly shows that the SSM is less robust to large forcing when the spectral quotient is low, resulting in a disassembling of the SSM.}
    \label{fig:manifold-transform}
\end{figure}

Figures~\ref{fig:ssm-robustness} and \ref{fig:manifold-transform} highlight the relationship between the spectral quotient and the robustness of the manifold under forcing. 
These figures show that the small-$\varepsilon$ assumption on the applied input discussed in Section~\ref{sec:SSMprelims} is nuanced and depends on the dynamics of the system -- specifically the spectral quotient. 
Notice that in Figure~\ref{fig:ssm-robustness}, increasing the spectral quotient allows us to increase the allowed forcing amplitude without significant change in error.
This would correspond to a \textit{preservation} of the periodic orbit in Figure~\ref{fig:manifold-transform} as we increase control effort to even larger amplitudes.

As shown in Figure~\ref{fig:manifold-transform}, for systems with low spectral quotient, the SSM quickly disassembles as we increase the amplitude of the input. 
We stress that while in this example the spectral quotient and damping are directly related, it is important to distinguish between the two when considering the previous discussion. 
Most structural dynamics applications feature small damping, but high spectral quotients because higher frequency modes exhibit higher damping ratios in comparison to low-frequency modes (see~\cite{Jain2018} for an analytic calculation of spectral quotients in a beam, for instance).  

The results presented here show promise for applying SSM-based control strategies to robotic systems with continuum-based models such as soft robots. 
In these systems the spectral quotients are expected to be high and are in fact, infinite in the continuum limit of structural finite-element models.
Indeed, in our recent work we show the applicability to higher-dimensional models with a data-driven SSM approach, as we discuss in \Cref{sub:future-works}.

\section{CONCLUSIONS AND FUTURE WORKS}
\label{sec:conclusions}

In this paper, we presented the first application of SSM theory for control.
In particular, we investigated the periodic setting and extended existing theoretical guarantees to synthesize optimal control policies for the purpose of periodic trajectory tracking. 
All existing ROM-based optimal control algorithms project the dynamics onto a linear subspace, resulting in the need to increase the dimension of the ROM to improve predictive capability for closed-loop control. 
In contrast, we reason directly about the nonlinearities during the reduction process using the powerful existence and uniqueness guarantees provided by the SSM.
We validated our approach on an illustrative example and provided insights on the persistence of the SSM under control inputs.


\subsection{Future Works}
\label{sub:future-works}
There are numerous extensions and applications of this work.
Direct applications to robotic platforms are appealing, such as highly-nonlinear soft robots or robotic fish with periodic tail actuation and/or undulation due to periodic muscle contraction \cite{wang2015fishaveraging}.
Similar to the pendulum example, we expect the ROM behavior to be coherent with the FOM dynamics even far away from the static equilibrium, enabling larger controllable workspaces, which are difficult to address with current piecewise-linear reduction techniques that have been investigated so far.
It would be interesting to apply more sophisticated control schemes which exploit the embodied intelligence of continuum robots and their in-resonant dynamics to produce hyper-efficient motions.

Our most recent work on applying data-driven Spectral Submanifold Reduction (SSMR) for nonlinear optimal control of a soft robot \cite{alora2022datadriven} shows the predictive capability of SSMs for real-world, high-dimensional robotic systems. 
By learning control-oriented models on low-dimensional SSMs, the proposed SSMR-based MPC approach outperforms both model-based and learning-based state-of-the-art methods in tracking performance and computational efficiency.
Further experimental and data-driven validation remains of great value to emphasize the applicability of the SSM-based approaches.
Open questions on generic time-dependent control inputs causing the SSM to lose its invariance provide future avenues of research.
For instance, characterizing model uncertainties can be useful for constraint-tightening schemes in safety-critical applications.

While data-driven SSMR shows significant promise, reliance on experimental data makes it difficult to apply these approaches in the design process. 
Thus, it remains worthwhile extending the model-based approach highlighted in this work for the design and control of high-dimensional, exotic robotic systems.






\printbibliography
\end{document}